\documentclass[12pt]{article}
\usepackage{amsmath,amstext, verbatim,remreset}
\usepackage[T2A]{fontenc}
\usepackage[cp1251]{inputenc}
\usepackage[english]{babel}
\usepackage{amsfonts,amssymb, amsmath}

\newcommand{\RR}{\mathbb R}
\newcommand{\NN}{\mathbb N}

\usepackage{mathrsfs}
\renewcommand{\leq}{\leqslant}
\renewcommand{\geq}{\geqslant}
\newtheorem{theorem}{Theorem}

\newtheorem{zam}{Remark}
\makeatletter
\renewcommand{\@begintheorem}[2]{\vskip 0.1 cm \noindent {\bf #1\ #2.\ }\begingroup\it}
\renewcommand{\@endtheorem}{\endgroup \vskip 0.1 cm}
\makeatother

\begin{document}

\begin{center}
{\bf \Large On dependence between norm of a function and norms of its derivatives of orders  $k$, $r-2$ and $r$, $0<k<r-2$.}

\bigskip

V. F. Babenko, O. V. Kovalenko

\bigskip

Oles Gonchar Dnipropetrovsk National University

\bigskip

\begin{abstract}
{ Necessary and sufficient conditions on the system of positive numbers $ M_{k_1}, M_{k_2}, M_{k_3}, M_{k_4}$, $0= k_1<k_2<k_3=r-2$, $k_4 = r$, which guarantee the existence of a function $x\in
L_{\infty,\infty}^r(\RR)$, such that $\|x^{(k_i)}\|_{\infty}=M_{k_i},\; i=1,2,3,4, $ are found.}
\end{abstract}
\end{center}

{\bf 1. Notations.  Statement of the problem. Known results.} Denote
by $L_\infty(\RR)$ the space of measurable essentially bounded functions $x\colon \RR\to \RR$ with norm
$$\|x\|=\|x\|_{L_\infty(\RR)}={\rm
ess\,sup}\left\{|x(t)|:t\in\RR\right\}.$$ For natural $r$ denote by $L_\infty^r(\RR)$ the space of functions $x\colon \RR\to
\RR$ such that the derivative $x^{(r-1)}, \; x^{(0)}=x,$ is locally absolute continuous and $x^{(r)}\in L_\infty(\RR)$. Let also
$L_{\infty,\infty}^r(\RR):=L_\infty^r(\RR)\bigcap L_\infty(\RR)$.

We will consider one of the cases of the following general problem stated by A.~N.~Kolmogorov \cite{Kol1}~--~\cite{Kol3} (in papers~\cite{Rodov1954}~--~\cite{Dzyadyk1974} it is mentioned, that this problem was stated by A.~N.~Kolmogorov in 1926).

{\bf Kolmogorov problem.} {\it Let a system of integers $0\leq k_1<k_2<...<k_d=r$ is given. Find necessary and sufficient conditions, on system of positive numbers $M_{k_1},...,M_{k_d}$, to guarantee the existence of a function $x\in L_{\infty,\infty}^r(\RR)$, such that
$$\left\|x^{(k_i)}\right\|=M_{k_i},\, i=1,...,d.$$}

Note (see remark in~\cite{Rodov1946}), that for case $d=2$, i. e. when we consider dependence between norm of the function ($k_1=0$) and norm of its $r$-th derivative, the solution of the problem is trivial: there exists a function that corresponds to any pair of positive numbers $M_{0}, \, M_{r}$.

In \cite{Kol1} -- \cite{Kol3} Kolmogorov solved this problem for $d = 3,\, k_1 = 0, \, 0<k_2<r$ (for $r=2$ and
$k=1$ this problem was solved by Adamar~\cite{had} earlier, for all cases with $r<5, k<r,$ and the case $r=5$ and $k=2$ --- by Shilov~\cite{shilov}). Kolmogorov showed that for positive numbers $M_0, M_k, M_r$, $0<k<r$, there exists a function $x \in
L_{\infty,\infty}^r(\RR)$ for which these numbers are the values of norms of the function, its $k$-th and its $r$-th derivative respectively if and only if  the following inequality holds 
\begin{equation}\label{-1}
M_k\leq \frac{\left\| \varphi _{r-k}\right\| }{\left\| \varphi _r\right\| ^{1-k/r}} M_0^{1-k/r}M_r^{k/r},
\end{equation}

\noindent where $\varphi _r$ --- is Euler perfect spline, i. e. $r$-th periodic integral with zero mean value on the period from the function
 $\varphi_0\left( t\right) ={\rm sgn}\sin t$. The solution of Kolmogorov problem for three numbers with arbitrary $k_1>0$ is contained, for example in~\cite[\S 9.1]{BBKP}.

 Solutions of Kolmogorov problem for $d>3$ are known in the following situations:
\begin{enumerate}
\item $k_1=0,\; k_2=r-2,\; k_3=r-1,\; k_4=r$ (Rodov~\cite{Rodov1946}).
\item $k_1=0<k_2<k_3=r-2,\; k_4=r-1, \;k_5=r$ (Rodov~\cite{Rodov1946}).
\item $k_1=0<k_2< k_3=r-1,\; k_4=r$ (Dzyadyk , Dubovik~\cite{Dzyadyk1975}).
\end{enumerate}

In \cite{Rodov1954} Rodov found sufficient conditions for the systems of positive numbers
$$M_0, M_1, M_2, M_5\;\; and\;\; M_0, M_1, M_2,M_3, M_4, M_5.$$

In the case of arbitrary $d>3$ the only known result belongs to Dzyadyk and Dubovik~\cite{Dzyadyk1975}.
They found sufficient conditions for existence of a functions $x \in L_{\infty,\infty}^r(\RR)$ with given values of norms of derivatives.

To prove the inequality~\eqref{-1} Kolmogorov proved statement, known as Kolmogorov comparison theorem. By comparison theorems one usually means statements that give estimation of some characteristics of a function $x(t)$ from some class using the corresponding characteristics of some fixed function. The last function can be counted as etalon or standard for given class; it is also called comparison function of given class.

Note, that in all mentioned partial solutions of Kolmogorov problem, the ideas, connected with comparison theorems, were essentially used. Kolmogorov comparison theorem itself and the method of its prove played important role for exact solutions of many extremal problems in approximation theory (see~\cite{korn1, korn2}).

The goal of this paper is to get the solution of Kolmogorov problem for the system of positive numbers $ M_{k_1}, M_{k_2}, M_{k_3}, M_{k_4}$,  $0=
k_1<k_2<k_3=r-2$, $k_4 = r$.

In the next paragraph we will introduce a family of splines and study their properties. In \S~3 we will prove the analogue of Kolmogorov comparison theorem for the case when norm of a function and norms of its derivatives of orders $r-2$ and $r$ are given. This theorem will be used not only for the solution of Kolmogorov problem, but has, in our opinion, independent interest. Here, as a corollary from comparison theorem, we will state Kolmogorov type inequality, that includes the norm of a function and the norms of its derivatives of orders $k$, $r-2$
and $r$. At last, in \S~4 the solution of Kolmogorov problem for the case $0= k_1<k_2<k_3=r-2$, $k_4 = r$ will be given.

\medskip

{\bf 2. Comparison functions and their properties.} Let $a\geq 0$. Define the function $\psi_1(a;t)$ in the following way. On interval $[0,a+2]$ set 
 $$
           \psi_1(a;t):=\left\{
           \begin{array}{rcl}
            t-1, &  t\in[0,1], \\
              0,   & t\in[1,a+1], \\
              t-a-1, & t\in [a+1,a+2].\\
           \end{array}
           \right.
  $$
  Continue function $\psi_1(a;t)$ evenly to the segment$[-(a+2),0]$, and then periodically with period $4+2a$ to the whole line.
  Note, that $\psi_1(a;t)\in L_{\infty,\infty}^1(\RR)$ and
  \begin{equation}\label{0.0}
  \left\|\psi'_1(a;\cdot)\right\|=1.
  \end{equation}

  For $r\in \NN$ denote by $\psi_r(a;t)$ $(r-1)$-th $(4+2a)$ -- periodic integral of the function $\psi_1(a;t)$ with zero mean value on period (so that, particularly, 
  $\psi'_r(a;t)=\psi_{r-1}(a;t)$).

  Note several properties of the function $\psi_r(a;t)$, which can be established either from definition, or analogously to corresponding properties of Euler perfect splines $\varphi_r$ (see, for example,~\cite[Chapter~5]{korn1},~\cite[Chapter~3] {korn2}). Fist of all note, that the function $\psi_2(a;t)$ is odd, $4+2a$--periodic,
  $$%\label{0}
           \psi_2(a;t):=\left\{
           \begin{array}{rcl}
            \frac 12 (t-1)^2- \frac 12, &  t\in[0,1], \\
              -\frac 12,   & t\in[1,a+1], \\
              \frac 12 (t-a-1)^2 -\frac 12, & t\in [a+1,a+2].\\
           \end{array}
           \right.
    $$
and
\begin{equation}\label{1'}
  \|\psi_2(a;\cdot)\|=\frac 12.
  \end{equation}
    Moreover, the functions $\psi_2(a;t)$ has exactly two zeroes on the period --- points $0$ and $a+2$.
    Hence the function $\psi_r(a;t)$ ($r\geq 2$) also has exactly two zeroes on the period: for all $k\in \NN$
    \begin{equation}\label{1}
    \psi_{2k}(a;0)=\psi_{2k}(a;a+2)=0,
    \end{equation}
        \begin{equation}\label{2}
     \psi_{2k+1}\left(a;1+\frac
     a2\right)=\psi_{2k+1}\left(a;3+\frac{3a}2\right)=0.
    \end{equation}

Hence, in turn, it follows that for $r\geq 3$ the function
$\psi_r(a;t)$ is strictly monotone between zeroes of its derivative, and the plot of the function $\psi_r(a;t)$ is convex on each interval of constant sign. Moreover, it is easy to see, that the plot of the function  $\psi_r(a;t)$ is symmetrical with respect to its zeroes and lines $t=t_0$, where $t_0$ is a zero  of $\psi'_r(a;t)$.
At last note, that if $\varphi_{\lambda, r}(t):=\lambda^{-r}\varphi_{r}(\lambda
 t)$ for $\lambda >0$, then $\psi_r(0;t)=\varphi_{\pi/2, r}(t)$.

 \medskip

 \begin{theorem}\label{th::2}
Let $r\in\NN$, $0<k<r-2$ and positive numbers $M_{k}$,
$M_{r-2}$, $M_r$ such, that Kolmogorov inequality
\begin{equation}\label{th2.0}
M_{r-2}\leq
\frac{\|\varphi_{2}\|}{\|\varphi_{r-k}\|^{\frac{2}{r-k}}}M_{k}^{\frac
2{r-k}}M_r^{\frac {r-k-2}{r-k}}
\end{equation}
holds are given.
Then there exist numbers $a,b,\lambda>0$ such  that for the function
$\Psi_{a,b,\lambda}(t):=b\psi(a;\lambda t)$ the following equalities hold
\begin{equation}\label{th2}
\left\|\Psi^{(s)}_{a,b,\lambda}\right\|=M_s,\,\,s\in\left\{k,r-2,r\right\}.
\end{equation}
In particular, for every function $x\in L^r_{\infty,\infty}(\RR)$
there exist numbers $a,b,\lambda>0$ such that
$$
\left\|\Psi^{(s)}_{a,b,\lambda}\right\|=\|
x^{(s)}\|,\,\,s\in\left\{k,r-2,r\right\}.
$$

\end{theorem}
\medskip

{\bf Proof.} Set
\begin{equation}\label{ab}
b:=\frac{M_r}{\lambda^r},\;\;\lambda:=\frac{\sqrt
{M_r}}{\sqrt{2M_{r-2}}}. \end{equation} Below in the proof of this theorem we count, that $b $ and $\lambda$ are chosen in the way that \eqref{ab} hold.
 Then in virtue of~\eqref{0.0}  and \eqref{1'} we get, that for all $a\geq 0$ we have
$$\left\|\Psi^{(r)}_{a,b,\lambda}\right\|=M_r,\;\;\;\left\|\Psi^{(r-2)}_{a,b,\lambda}\right\|=M_{r-2}.$$

It is clear that for all $k=1,2,...,r-3$ the function
$\left\|\Psi^{(k)}_{a,b,\lambda}\right\|$ continuously depends on $a\in[0,\infty)$, increases on this interval and 
$$\lim\limits_{a\to+\infty}\left\|\Psi^{(k)}_{a,b,\lambda}\right\|=\infty.$$

Since the inequality \eqref{th2.0} turns into equality for the function 
$\Psi^{(k)}_{0,b,\lambda}$, we have
$$\left\|\Psi^{(k)}_{0,b,\lambda}\right\|=\left(\frac{\|\varphi_{r-k}\|^{\frac
2{r-k}}}{\|\varphi_2\|}\cdot
\frac{M_{r-2}}{M_r^{\frac{r-k-2}{r-k}}}\right)^{\frac{r-k}2}\leq
M_k.$$ Hence there exists $a\geq 0$ such that
$$
\left\|\Psi^{(k)}_{a,b,\lambda}\right\|=M_k.
$$

 The theorem is proved.

\medskip

\begin{zam}\label{z::1}
For numbers $M_k,M_{r-2},M_r$, satisfying the inequality~\eqref{th2.0} by $\Psi_r(M_{k},M_{r-2},M_{r};t)$ we will denote  the function $\Psi_{a,b,\lambda}(t)$ from theorem~$1$, with parameters $a, b, \lambda$ chosen in the way that equalities
\eqref{th2} hold.
\end{zam}

{\bf 3. Comparison theorem and Kolmogorov inequality analogue.}
The following theorem is an analogue of Kolmogorov comparison theorem in the case when norms of a function and its derivatives of orders $(r-2)$ and $r$ are given.

 \begin{theorem}\label{th::0}
 Let $r\in\NN$, $0=k_1<k_2<k_3=r-2$, $k_4=r$ and $x\in L_{\infty,\infty}^r(\RR)$ are given.
  Let the numbers $a,b,\lambda>0$ are such, that for a function $\Psi_{a,b,\lambda}(t)$ the following equalities hold
 \begin{equation}\label{4'}
\left\|x^{(k_i)}\right\|\le
\left\|\Psi^{(k_i)}_{a,b,\lambda}\right\|,\,\,i=1,3,4.
 \end{equation}
 If points $\tau$ and $\xi$ are such that $x(\tau) = \Psi_{a,b,\lambda}(\xi)$, then
\begin{equation}\label{6'}\left|x'(\tau)\right|\le\left|\Psi'_{a,b,\lambda}(\xi)\right|.
 \end{equation}
 \end{theorem}

{\bf Proof.}  For brevity we will write $\Psi (t)$ instead of
$\Psi_{a,b,\lambda}(t)$ in the proof of this theorem. Considering, if necessary, the function  $-x(t)$ instead of $x(t)$ and function $-\Psi (t)$ instead of 
$\Psi (t)$, we can count that $x'(\tau)>0$ and
 \begin{equation}\label{7}
 \Psi'(\tau)>0.
 \end{equation}
Moreover, considering appropriate shift $\Psi (\cdot +\alpha)$ of the function $\Psi$, we can count that $\tau = \xi$, i.~e.
 \begin{equation}\label{6}
 x(\tau)=\Psi(\tau).
 \end{equation}
Assume, that~\eqref{6} holds, but instead of the inequality~(\ref{6'}) (with $\xi =\tau$) the inequality
$$\left|x'(\tau)\right|>\left|\Psi'(\tau)\right| $$ holds.

Denote by $(\tau_1,\tau_2)$ the smallest interval which contains  $\tau$ on which the function $\Psi$ is monotone and such that $\Psi'(\tau_1)=\Psi'(\tau_2)=0$.
In virtue of the assumption there exists a number $\delta>0$ such that $x'(t)>\Psi'(t)$ for all $t\in (\tau-\delta,\tau+\delta)$, and hence in virtue of~\eqref{6} $x(\tau+\delta)>\Psi(\tau+\delta)$ and $x(\tau-\delta)<\Psi(\tau-\delta)$.

Choose $\varepsilon>0$ so small, that for a function
$x_\varepsilon(t):=(1-\varepsilon)x(t)$ the following inequalities hold:
$x_\varepsilon(\tau+\delta)>\Psi(\tau+\delta)$ and
$x_\varepsilon(\tau-\delta)<\Psi(\tau-\delta)$. In virtue of the conditions~\eqref{4'} and~\eqref{7} we have
$$x_\varepsilon(\tau_1)>\Psi(\tau_1),\;\; x_\varepsilon(\tau_2)<\Psi(\tau_2).$$
Hence on the interval $(\tau_1,\tau_2)$ the difference
$\Delta_\varepsilon(t):=x_\varepsilon(t)-\Psi(t)$ has at least 3 sign changes.

It is easy to see, that there exists a sequence of functions $\mu_N\in
C^\infty({\mathbb R})$, $N\in {\mathbb N}$ with the following properties:
\begin{enumerate}
\item $\mu_N(t)=1$ on interval $[\tau_1,\tau_2]$; $\|\mu_N\|=1$;
\item $\mu_N(t)=0$ for all $t$ outside the interval
$[\tau_1-N\cdot\frac{2+a}{\lambda};\tau_1+N\cdot\frac{2+a}{\lambda}]$;
\item for all
$k=1,2,\dots,r$
$$\max\limits_{j=\overline{1,k}}\left\|\mu_N^{(j)}\right\|<
\varepsilon \|x_\varepsilon^{(k)}\|\left(\sum\limits_{i=1}^k
C_k^i\left\|x_\varepsilon^{(k-i)}\right\|\right)^{-1},$$ if $N$
is enough big.
\end{enumerate}

Below we count that $N$ is chosen enough big, so that the property 3 holds.

Set $$x_N(t):=x_\varepsilon(t)\cdot\mu_N(t),$$
and
$$\Delta_N(t):=\Psi(t)-x_N(t).$$ Then
$$x_N(t)=x_\varepsilon(t),\;if \; t\in[\tau_1,\tau_2],$$
\begin{equation}\label{7'}
\Delta_N(t)=\Psi(t),if \,|t-\tau_1|\geq N\cdot\frac{a+2}\lambda
\end{equation}
and
$$
\|x_N\|\leq\|x_\varepsilon\|=(1-\varepsilon)\|x\|\leq(1-\varepsilon)\|\Psi\|.
$$

Moreover, for  $k=1,\dots,r$ $$\left|x_N^{(k)}(t)\right|=
\left|\left[x_\varepsilon(t)\mu_N(t)\right]^{(k)}\right|=
\left|\sum\limits_{i=0}^k C_k^i
x_\varepsilon^{(k-i)}(t)\mu_N^{(i)}(t)\right|\leq
$$
$$
\leq \left\|x_\varepsilon^{(k)}\right\|+\sum\limits_{i=1}^k
C_k^i\left\|x_\varepsilon^{(k-i)}\right\|\|\mu_N^{(i)}\|.$$ Hence,
in virtue of property 3 of the function $\mu_N$ and the choice of the number $N$, we get
$$
\left\|x_N^{(k)}\right\|<
\left\|x_\varepsilon^{(k)}\right\|+\varepsilon\left\|x^{(k)}_\varepsilon\right\|=(1-\varepsilon)\left\|x^{(k)}\right\|+\varepsilon\left\|x^{(k)}\right\|=\left\|x^{(k)}\right\|.
$$

For $t\in [\tau_1,\tau_2]$ we have
$\Delta_N=\Psi(t)-x_\varepsilon(t),$ and hence the function $\Delta_N(t)$
has at least three sign changes on the interval $[\tau_1,\tau_2]$. At each of the rest monotonicity intervals of the function $\Psi$ the function 
$\Delta_N$ has at least one sign change. Hence on the interval
$\left[\tau_1-N\cdot\frac{a+2}{\lambda},\tau_1+N\cdot\frac{a+2}{\lambda}\right]$
the function $\Delta_N(t)$ has at least $2N+2$ sign changes. Moreover, in virtue of~\eqref{1}, \eqref{2} and \eqref{7'} for all
$i=1,2,\dots,\left[\frac{r-1}{2}\right]$ the following equalities hold
\begin{equation}\label{10}
\Delta_N^{(2i-1)}\left(\tau_1-N\cdot\frac{a+2}\lambda\right)=\Delta_N^{(2i-1)}
\left(\tau_1+N\cdot\frac{a+2}\lambda\right)=0.
\end{equation}

Applying Rolle's theorem and counting~\eqref{10} we have that the function $\Delta_N^{(r-2)}(t)$ has at least $2N+2$ zeroes on the interval
$\left[\tau_1-N\cdot\frac{a+2}{\lambda},\tau_1+N\cdot\frac{a+2}{\lambda}\right]$.
Hence on some monotonicity interval
$[\alpha,\alpha+\frac{2+a}\lambda]\subset
\left[\tau_1-N\cdot\frac{a+2}{\lambda},
\tau_1+N\cdot\frac{a+2}{\lambda}\right]$ ($\alpha:= \frac
1\lambda+k\cdot\frac{a+2}\lambda,\,k\in \NN$) of the function
$\Psi^{(r-2)}(t)=b\lambda^{-2}\psi (a, \lambda t)$ the function
$\Delta_N^{(r-2)}(t)$ changes sign at least three times. But then the difference 
$$
b\lambda^{-2}\psi (0, \lambda t)-x^{(r-2)}_N(t)
$$
changes the sign at least three times on some monotonicity interval of the function $b\lambda^{-2}\psi (0,
\lambda t)$ too. However this contradicts to the Kolmogorov comparison theorem (see, for example,~\cite[Statement 5.5.3]{korn2}) because the Euler spline $b\lambda^{-2}\psi (0, \lambda
t)$ is comparison function for the function $x^{(r-2)}_N(t)$.

As a corollary of the theorem~\ref{th::0} we get the following theorem, which can be viewed as Kolmogorov type inequality, which estimates the norm of $k$ -- th derivative of a function, by the norms of the function and its derivatives of orders $(r-2)$ and $r$.

\begin{theorem}\label{th::1}
 Let $r\in\NN$, $0=k_1<k_2<k_3=r-2$, $k_4=r$ and $x\in L_{\infty,\infty}^r(\RR)$ are given.
 Let numbers $a,b,\lambda>0$ are such that for a function $\Psi_{a,b,\lambda}(t)$ the inequalities~\eqref{4'} hold.
 Then $$\left\|x^{(k_2)}\right\|\leq\left\|\Psi^{(k_2)}_{a,b,\lambda}\right\|.$$
 \end{theorem}

\medskip

{\bf 4. Solution of Kolmogorov problem for the case when $0=k_1<k_2<k_3=r-2$, $k_4=r$}.
\begin{theorem}
Let integers $r\geq 4$,  $0=k_1<k_2<k_3=r-2$, $k_4=r$ and real numbers $M_{k_1},M_{k_2},M_{k_3},M_{k_4}>0$ are given. There exists a function $x\in L_{\infty,\infty}^r(\RR)$, such that 
\begin{equation}\label{11}
\left\|x^{(k_i)}\right\|=M_{k_i},\,\,i=1,2,3,4
\end{equation} if and only if the following inequalities hold
$$a)\,\,M_{r-2}\leq \frac{\|\varphi_{2}\|}{\|\varphi_{r-k}\|^{\frac{2}{r-k}}}M_{k}^{\frac 2{r-k}}M_r^{\frac {r-k-2}{r-k}},$$
$$b)\,\, M_0\geq \|\Psi_r(M_{k_2},M_{r-2},M_r)\|,\phantom{aaaa}$$ where the function $\Psi_r$ is defined in remark~$\ref{z::1}$.
\end{theorem}

The necessity of the condition~$a)$ follows from Kolmogorov inequality, the necessity of the condition~$b)$ follows from theorem~\ref{th::1}.

To prove the sufficiency it is enough to note that in the case when conditions $a)$ and $b)$ hold for a function
$$x(t):=\Psi_r(M_{k_2},M_{r-2},M_r;t)+M_0-\left\|\Psi_r(M_{k_2},M_{r-2},M_r)\right\|$$ the equalities~\eqref{11} hold.

Theorem is proved.

\end{document}

\leftline{УДК 517.5}

\begin{center}
{\bf \Large О зависимости между нормой функции и нормами ее
производных порядка $k$, $r-2$ и $r$, $0<k<r-2$.}
\end{center}

\begin{center}
{В. Ф. Бабенко, О. В. Коваленко }

Днепропетровский национальный университет им. О. Гончара

\end{center}
\begin{abstract}
{ Найдены необходимые и достаточные условия на систему положительных
чисел $ M_{k_1}, M_{k_2}, M_{k_3}, M_{k_4}$,  $0= k_1<k_2<k_3=r-2$,
$k_4 = r$, гарантирующие существование функции $x\in
L_{\infty,\infty}^r(\RR)$, такой, что $
\|x^{(k_i)}\|_{\infty}=M_{k_i},\; i=1,2,3,4. $}

\end{abstract}

{\bf 1. Обозначения. Постановка задачи. Известные результаты.} Через
$L_\infty(\RR)$ будем обозначать пространство измеримых и
сущест\-вен\-но ограничен\-ных функций $x\colon \RR\to \RR$ с нормой
$$\|x\|=\|x\|_{L_\infty(\RR)}={\rm
ess\,sup}\left\{|x(t)|:t\in\RR\right\}.$$ Для натурального $r$ через
$L_\infty^r(\RR)$ обозначим пространство  функций $x\colon \RR\to
\RR$ таких, что производная $x^{(r-1)}, \; x^{(0)}=x,$ локально
абсолютно непрерывна, и $x^{(r)}\in L_\infty(\RR)$. Пусть также
$L_{\infty,\infty}^r(\RR):=L_\infty^r(\RR)\bigcap L_\infty(\RR)$.

Мы будем рассматривать один из случаев следующей общей задачи, постав\-лен\-ной А. Н. Колмогоровым \cite{Kol1} -- \cite{Kol3} (в работах \cite{Rodov1954} -- \cite{Dzyadyk1974} отмечено, что эта задача поставлена А.Н. Колмогоровым в 1926 г.).

{\bf Задача Колмогорова.} {\it Пусть задана  система целых чисел
$0\leq k_1<k_2<...<k_d=r$. Найти необходимые и достаточные условия,
которым должна удовлетворять  систе\-ма положительных чисел
$M_{k_1},...,M_{k_d}$, для того, чтобы существовала функция $x\in
L_{\infty,\infty}^r(\RR)$, такая, что
$$\left\|x^{(k_i)}\right\|=M_{k_i},\, i=1,...,d.$$}

Отметим (см. замечание в ~\cite{Rodov1946}), что для $d=2$, т. е.
когда речь идет о зависимости между нормой функции ($k_1=0$) и
нормой  ее $r$-й производной, решение задачи тривиально: существуют
функции, соответствующие любой паре положительных чисел $M_{0}, \, M_{r}$.

В \cite{Kol1} -- \cite{Kol3} Колмогоров привел формулировку и
решение этой задачи для $d = 3,\, k_1 = 0, \, 0<k_2<r$ (для $r=2$ и
$k=1$ эта задача была решена ранее Ж. Адамаром~\cite{had}, а для
всех случаев $r<5, k<r,$ и случая $r=5$ и $k=2$ -- Г.Е.
Шиловым~\cite{shilov}). Колмогоров показал, что для положительных
чисел $M_0, M_k, M_r$, $0<k<r$, существует функция $x \in
L_{\infty,\infty}^r(\RR)$ для которой эти три числа являются
значениями норм функции, её $k$-ой и её $r$-ой производной
соответственно тогда, и только тогда, когда выполняется неравенство
\begin{equation}\label{-1}
M_k\leq \frac{\left\| \varphi _{r-k}\right\| }{\left\| \varphi _r\right\| ^{1-k/r}} M_0^{1-k/r}M_r^{k/r},
\end{equation}
где
 $\varphi _r$ -- это эйлеров идеальный сплайн, т. е. $r$-ый периодический интеграл с нулевым средним значением на периоде
от функции $\varphi_0\left( t\right) ={\rm sgn}\sin t$. Решение
задачи Колмогорова для трёх чисел с произвольным $k_1>0$ содержится,
например, в ~\cite[\S 9.1]{BBKP}.

 Известные результаты по решению задачи Колмогорова при $d>3$ таковы:
\begin{enumerate}
\item $k_1=0,\; k_2=r-2,\; k_3=r-1,\; k_4=r$ (Родов ~\cite{Rodov1946}).
\item $k_1=0<k_2<k_3=r-2,\; k_4=r-1, \;k_5=r$ (Родов ~\cite{Rodov1946}).
\item $k_1=0<k_2< k_3=r-1,\; k_4=r$ (Дзядык, Дубовик ~\cite{Dzyadyk1975}).
\end{enumerate}

В \cite{Rodov1954} Родов нашёл достаточные условия для систем
положительных чисел
$$M_0, M_1, M_2, M_5\;\; {\rm \textcyrillic{и}}\;\; M_0, M_1, M_2,M_3, M_4, M_5.$$

 В случае произвольного $d>3$ единственный известный результат принад\-лежит
 Дзядыку и Дубовику~\cite{Dzyadyk1975}.
Ими получены достаточные условия для существования функции $x \in
L_{\infty,\infty}^r(\RR)$ с заданными значениями норм
произ\-вод\-ных фиксированных порядков.

Для доказательства неравенства \eqref{-1} Колмогоровым было доказано
утверж\-дение, известное как теорема сравнения Колмогорова.
Теоремами сравнения называют утверждения, которые дают оценку той
или иной характеристики функции $x(t)$ из некоторого класса через
соответствующую характеристику некоторой фиксированной функции.
Последнюю функцию можно считать эталонной или стандартной для для
данного класса; её также называют функ\-цией сравнения для данного
класса.

Отметим, что во всех приведённых частных решениях задачи Колмогорова
существенно использовались идеи, связанные с теоремами сравнения.
Как сама теорема сравнения Колмогорова, так и метод её
доказательства сыграли большую роль при точном решении многих
экстремальных задач теории приближений (см.~\cite{korn1, korn2}).

Цель данной статьи -- решение задачи Колмогорова для системы
положи\-тельных чисел $ M_{k_1}, M_{k_2}, M_{k_3}, M_{k_4}$,  $0=
k_1<k_2<k_3=r-2$, $k_4 = r$.

В следующем параграфе мы введем некоторое семейство сплайнов и
изу\-чим их свойства.  В параграфе 3 мы докажем  аналог теоремы
сравнения Колмогорова для случая, когда заданы норма функции и нормы
ее произ\-вод\-ных порядков $r-2$ и $r$. Эта теорема не только будет
использована для решеня задачи Колмогорова, но и представля\-ет, на
наш взгляд, самостоятельный интерес. Здесь же в качестве следствия
из теоремы сравнения мы приведем неравенство типа Колмогорова,
включаю\-щее норму функции и нормы ее производных порядка $k$, $r-2$
и $r$. Наконец, в параграфе 4 будет приведено решение задачи
Колмогорова для случая $0= k_1<k_2<k_3=r-2$, $k_4 = r$.

\medskip

{\bf 2. Функции сравнения и их свойства.} Пусть $a\geq 0$. Определим
функ\-цию $\psi_1(a;t)$ следующим образом. На отрезке $[0,a+2]$
положим

 $$
           \psi_1(a;t):=\left\{
           \begin{array}{rcl}
            t-1, &  t\in[0,1], \\
              0,   & t\in[1,a+1], \\
              t-a-1, & t\in [a+1,a+2].\\
           \end{array}
           \right.
  $$
  Продолжим функцию $\psi_1(a;t)$ чётным образом на отрезок $[-(a+2),0]$, а затем
  периодически с периодом $4+2a$ на всю ось.
  Заметим, что $\psi_1(a;t)\in L_{\infty,\infty}^1(\RR)$ и
  \begin{equation}\label{0.0}
  \left\|\psi'_1(a;\cdot)\right\|=1.
  \end{equation}

  Для $r\in \NN$ обозначим через $\psi_r(a;t)$ $(r-1)$-ю $(4+2a)$ -- периодическую
  первообразную функции $\psi_1(a;t)$ с нулевым средним на периоде (так что, в частности,
  $\psi'_r(a;t)=\psi_{r-1}(a;t)$).

  Приведем некоторые свойства функции $\psi_r(a;t)$, которые
  нетрудно устано\-вить либо непосредственно, либо по аналогии со
  свойствами эйлеровых иде\-альных сплайнов $\varphi_r$ (см., например,~\cite[гл. 5]{korn1},~\cite[гл. 3]
  {korn2}). Отметим, прежде всего, что функ\-ция $\psi_2(a;t)$ -- нечётна, имеет период $4+2a$,
  $$%\label{0}
           \psi_2(a;t):=\left\{
           \begin{array}{rcl}
            \frac 12 (t-1)^2- \frac 12, &  t\in[0,1], \\
              -\frac 12,   & t\in[1,a+1], \\
              \frac 12 (t-a-1)^2 -\frac 12, & t\in [a+1,a+2].\\
           \end{array}
           \right.
    $$
и
\begin{equation}\label{1'}
  \|\psi_2(a;\cdot)\|=\frac 12.
  \end{equation}
    Кроме того функция $\psi_2(a;t)$ имеет ровно два нуля на периоде -- точки $0$ и $a+2$.
    Следовательно функция   $\psi_r(a;t)$ при $r\geq 2$ также имеет ровно два нуля на
    периоде: для любого $k\in \NN$
    \begin{equation}\label{1}
    \psi_{2k}(a;0)=\psi_{2k}(a;a+2)=0,
    \end{equation}
        \begin{equation}\label{2}
     \psi_{2k+1}\left(a;1+\frac
     a2\right)=\psi_{2k+1}\left(a;3+\frac{3a}2\right)=0.
    \end{equation}

%    Ряд Фурье функции $\psi_1(a;t)$ имеет вид $$\psi_1(a;t)=\sum_{k=1}^{\infty}a_k\cos\frac{2\pi k t}{2a+4},$$
%где $$a_k=\frac 2{a+2}\int\limits_0^{a+2}\psi_1(a;t)\cos\frac{2\pi k t}{2a+4}dt=\frac{2}{a+2}\left\{\int\limits_0^1 (t-1)\cos\frac{2\pi k t}{2a+4}dt+\right.$$
%
%$$+\left.\int\limits_{a+1}^{a+2}(t-a-1)\cos\frac{2\pi k t}{2a+4}dt\right\}=-\frac{2}{\pi k}\left\{\int\limits_0^1 \sin\frac{\pi k t}{a+2}dt+\int\limits_{a+1}^{a+2} \sin\frac{\pi k t}{a+2}dt\right\}=$$
%
%$$=\frac{2(a+2)}{(\pi k)^2}\left\{\cos\frac{\pi k}{a+2} -1 + \cos\pi k - \cos\frac{\pi (a+1) k}{a+2}\right\}=$$ $$=\frac{2(a+2)}{(\pi k)^2}\left\{2\sin\frac{\pi k}{2}\sin\frac{a\pi k}{2(a+2)} + (-1)^k-1\right\}.$$
%
% Если $k$ четное, то $a_k=0$. Поэтому $$\psi_1(a;t)=\frac{4(a+2)}{\pi^2}\sum_{k=0}^{\infty}\frac{(-1)^k\sin\frac{a\pi (2k+1)}{2(a+2)}-1}{(2k+1)^2}\cos\frac{\pi (2k+1) t}{a+2}.$$
%
% Отсюда интегрируя мы получаем, что
% \begin{equation}\label{3}
% \psi_r(a;t)=(-1)^\frac{r-1}2\frac{4(a+2)^r}{\pi^{r+1}}\sum_{k=0}^{\infty}\frac{(-1)^k\sin\frac{a\pi (2k+1)}{2(a+2)}-1}{(2k+1)^{r+1}}\cos\frac{\pi (2k+1) t}{a+2},\,\,r=2i+1,\,i=0,1,...
% \end{equation}
%и
% \begin{equation}\label{4}
% \psi_r(a;t)=(-1)^\frac{r+2}2\frac{4(a+2)^r}{\pi^{r+1}}\sum_{k=0}^{\infty}\frac{(-1)^k\sin\frac{a\pi (2k+1)}{2(a+2)}-1}{(2k+1)^{r+1}}\sin\frac{\pi (2k+1) t}{a+2},\,\,r=2i,\,i=1,2,...
% \end{equation}
%
% Из равенств \eqref{3} и \eqref{4} следует справедливость \eqref{1} и \eqref{2}.

Отсюда, в свою очередь, следует, что при $r\geq 3$ функция
$\psi_r(a;t)$ строго монотонна между нулями своей
 производной,
 а график функции $\psi_r(a;t)$ является выпуклым на каждом промежутке
 знакопостоянства. Кроме того, как легко видеть, график
 $\psi_r(a;t)$ симметричен относительно ее нулей, а также
 относительно прямых вида $t=t_0$, где $t_0$ -- нуль $\psi'_r(a;t)$.
 Наконец отметим, что если $\varphi_{\lambda, r}(t):=\lambda^{-r}\varphi_{r}(\lambda
 t)$ для $\lambda >0$, то $\psi_r(0;t)=\varphi_{\pi/2, r}(t)$.

 \medskip

 \begin{theorem}\label{th::2}
Пусть $r\in\NN$, $0<k<r-2$ и заданы положительные числа $M_{k}$,
$M_{r-2}$, $M_r$ такие,
 что выполняется неравенство Колмогорова
\begin{equation}\label{th2.0}
M_{r-2}\leq
\frac{\|\varphi_{2}\|}{\|\varphi_{r-k}\|^{\frac{2}{r-k}}}M_{k}^{\frac
2{r-k}}M_r^{\frac {r-k-2}{r-k}}.
\end{equation}
Тогда существуют числа $a,b,\lambda>0$ такие, что для функции
$\Psi_{a,b,\lambda}(t):=b\psi(a;\lambda t)$ выполняются равенства
\begin{equation}\label{th2}
\left\|\Psi^{(s)}_{a,b,\lambda}\right\|=M_s,\,\,s\in\left\{k,r-2,r\right\}.
\end{equation}
В частности, для любой функции $x\in L^r_{\infty,\infty}(\RR)$
существуют числа $a,b,\lambda>0$ такие, что
$$
\left\|\Psi^{(s)}_{a,b,\lambda}\right\|=\|
x^{(s)}\|,\,\,s\in\left\{k,r-2,r\right\}.
$$

\end{theorem}
\medskip

{\bf Доказательство.} Положим
\begin{equation}\label{ab}
b:=\frac{M_r}{\lambda^r},\;\;\lambda:=\frac{\sqrt
{M_r}}{\sqrt{2M_{r-2}}}. \end{equation} Везде ниже в ходе данного
доказательства считаем, что $b $ и $\lambda$ выбраны в соответствии
 с \eqref{ab}.
 Тогда в силу \eqref{0.0}  и \eqref{1'}
получаем, что при всех $a\geq 0$ будет
$$\left\|\Psi^{(r)}_{a,b,\lambda}\right\|=M_r,\;\;\;\left\|\Psi^{(r-2)}_{a,b,\lambda}\right\|=M_{r-2}.$$

Ясно, что для всех $k=1,2,...,r-3$ функция
$\left\|\Psi^{(k)}_{a,b,\lambda}\right\|$ непрерывно зависит от
$a\in[0,\infty)$, возрастает на этом промежутке, причем
$$\lim\limits_{a\to+\infty}\left\|\Psi^{(k)}_{a,b,\lambda}\right\|=\infty.$$

Поскольку  неравенство
 \eqref{th2.0} обращается в равенство для функции
$\Psi^{(k)}_{0,b,\lambda}$, получаем
$$\left\|\Psi^{(k)}_{0,b,\lambda}\right\|=\left(\frac{\|\varphi_{r-k}\|^{\frac
2{r-k}}}{\|\varphi_2\|}\cdot
\frac{M_{r-2}}{M_r^{\frac{r-k-2}{r-k}}}\right)^{\frac{r-k}2}\leq
M_k.$$ Следовательно, существует $a\geq 0$ такое, что
$$
\left\|\Psi^{(k)}_{a,b,\lambda}\right\|=M_k.
$$

 Теорема доказана.

\medskip

\begin{zam}\label{z::1}
Для чисел $M_k,M_{r-2},M_r$, удовлетворяющих неравенству
\eqref{th2.0} будем обозначать через $\Psi_r(M_{k},M_{r-2},M_{r};t)$
функцию $\Psi_{a,b,\lambda}(t)$ из теоремы~$1$, параметры $a, b,
\lambda$ которой выбраны так, что выполняются соотношения
\eqref{th2}.
\end{zam}

{\bf 3. Теорема сравнения и аналог неравенства Колмогорова.}
Следу\-ющая теорема является аналогом теоремы сравнения Колмогорова
в случае, когда заданы норма функции, ее $(r-2)$ - й и $r$ - й
производных.

 \begin{theorem}\label{th::0}
 Пусть $r\in\NN$ и $0=k_1<k_2<k_3=r-2$, $k_4=r$ и $x\in L_{\infty,\infty}^r(\RR)$.
  Пусть числа $a,b,\lambda>0$ таковы, что для функции $\Psi_{a,b,\lambda}(t)$
  выпол\-няются соотношения
 \begin{equation}\label{4'}
\left\|x^{(k_i)}\right\|\le
\left\|\Psi^{(k_i)}_{a,b,\lambda}\right\|,\,\,i=1,3,4.
 \end{equation}
 Если точки $\tau$ и $\xi$ таковы, что $x(\tau) = \Psi_{a,b,\lambda}(\xi)$, то
\begin{equation}\label{6'}\left|x'(\tau)\right|\le\left|\Psi'_{a,b,\lambda}(\xi)\right|.
 \end{equation}
 \end{theorem}

{\bf Доказательство.} Для сокращения записей в ходе данного
доказательства мы будем писать $\Psi (t)$ вместо
$\Psi_{a,b,\lambda}(t)$. Рассматривая при необходимости
функ\-цию $-x(t)$ вместо функции $x(t)$ и функцию $-\Psi (t)$ вместо
$\Psi (t)$, мы можем считать, что $x'(\tau)>0$ и
 \begin{equation}\label{7}
 \Psi'(\tau)>0.
 \end{equation}
Кроме того, рассматривая подходящий сдвиг $\Psi (\cdot +\alpha)$
функции $\Psi$ можем считать, что  и $\tau = \xi$, т.е.
 \begin{equation}\label{6}
 x(\tau)=\Psi(\tau).
 \end{equation}

Предположим, что~(\ref{6}) имеет место, но при этом вместо
неравенства~(\ref{6'}) (с $\xi =\tau$) имеет место неравенство
$$\left|x'(\tau)\right|>\left|\Psi'(\tau)\right|.
 $$

 Обозначим через $(\tau_1,\tau_2)$ наименьший интервал монотонности $\Psi$, содержа\-щий
 точку $\tau$ и такой, что $\Psi'(\tau_1)=\Psi'(\tau_2)=0$.
В силу сделанного предположе\-ния существует число $\delta>0$ такое,
что $x'(t)>\Psi'(t)$ для всех $t\in (\tau-\delta,\tau+\delta)$, а
значит в силу \eqref{6} $x(\tau+\delta)>\Psi(\tau+\delta)$ и
$x(\tau-\delta)<\Psi(\tau-\delta)$.

Выберем $\varepsilon>0$ настолько малым, чтобы для функции
$x_\varepsilon(t):=(1-\varepsilon)x(t)$ выполнялись неравенства
$x_\varepsilon(\tau+\delta)>\Psi(\tau+\delta)$ и
$x_\varepsilon(\tau-\delta)<\Psi(\tau-\delta)$. В силу \eqref{4'} и
\eqref{7} будет
$$x_\varepsilon(\tau_1)>\Psi(\tau_1),\;\; x_\varepsilon(\tau_2)<\Psi(\tau_2).$$
Таким образом на промежутке $(\tau_1,\tau_2)$ разность
$\Delta_\varepsilon(t):=x_\varepsilon(t)-\Psi(t)$ будет иметь не
менее трёх перемен знака.

Как легко видеть,  существует проследовательность функций $\mu_N\in
C^\infty(\RR)$, $N\in \NN$, со следующими свойствами:
\begin{enumerate}
\item $\mu_N(t)=1$ на промежутке $[\tau_1,\tau_2]$; $\|\mu_N\|=1$;
\item $\mu_N(t)=0$ для всех $t$ вне промежутка
$[\tau_1-N\cdot\frac{2+a}{\lambda};\tau_1+N\cdot\frac{2+a}{\lambda}]$;
\item для всех
$k=1,2,\dots,r$
$$\max\limits_{j=\overline{1,k}}\left\|\mu_N^{(j)}\right\|<
\varepsilon \|x_\varepsilon^{(k)}\|\left(\sum\limits_{i=1}^k
C_k^i\left\|x_\varepsilon^{(k-i)}\right\|\right)^{-1},$$ если $N$
достаточно велико.
\end{enumerate}

Ниже считаем, что $N$ выбрано настолько большим, что свойство 3
выполнено.

Положим $$x_N(t):=x_\varepsilon(t)\cdot\mu_N(t),$$
и
$$\Delta_N(t):=\Psi(t)-x_N(t).$$ Тогда
$$x_N(t)=x_\varepsilon(t),\;\textcyrillic{если }\; t\in[\tau_1,\tau_2],$$
\begin{equation}\label{7'}
\Delta_N(t)=\Psi(t),\textcyrillic{если }\,|t-\tau_1|\geq N\cdot\frac{a+2}\lambda
\end{equation}
и
$$
\|x_N\|\leq\|x_\varepsilon\|=(1-\varepsilon)\|x\|\leq(1-\varepsilon)\|\Psi\|.
$$

Кроме того для  $k=1,\dots,r$ $$\left|x_N^{(k)}(t)\right|=
\left|\left[x_\varepsilon(t)\mu_N(t)\right]^{(k)}\right|=
\left|\sum\limits_{i=0}^k C_k^i
x_\varepsilon^{(k-i)}(t)\mu_N^{(i)}(t)\right|\leq
$$
$$
\leq \left\|x_\varepsilon^{(k)}\right\|+\sum\limits_{i=1}^k
C_k^i\left\|x_\varepsilon^{(k-i)}\right\|\|\mu_N^{(i)}\|.$$ Отсюда,
с учетом свойства 3 функции $\mu_N$ и выбора числа $N$, получаем
$$
\left\|x_N^{(k)}\right\|<
\left\|x_\varepsilon^{(k)}\right\|+\varepsilon\left\|x^{(k)}_\varepsilon\right\|=(1-\varepsilon)\left\|x^{(k)}\right\|+\varepsilon\left\|x^{(k)}\right\|=\left\|x^{(k)}\right\|.
$$

Для $t\in [\tau_1,\tau_2]$ имеем
$\Delta_N=\Psi(t)-x_\varepsilon(t),$ а значит функция $\Delta_N(t)$
имеет не менее  трех перемен знака на отрезке $[\tau_1,\tau_2]$. На
каждом из остальных промежутков монотонности функции $\Psi$ функция
$\Delta_N$ имеет не менее одной перемены знака. Таким образом на
промежутке
$\left[\tau_1-N\cdot\frac{a+2}{\lambda},\tau_1+N\cdot\frac{a+2}{\lambda}\right]$
функция $\Delta_N(t)$ имеет не менее $2N+2$ перемен знака. Кроме
того, в силу \eqref{1},\eqref{2} и \eqref{7'} для всех
$i=1,2,\dots,\left[\frac{r-1}{2}\right]$ справедливы равенства
\begin{equation}\label{10}
\Delta_N^{(2i-1)}\left(\tau_1-N\cdot\frac{a+2}\lambda\right)=\Delta_N^{(2i-1)}
\left(\tau_1+N\cdot\frac{a+2}\lambda\right)=0.
\end{equation}

Применяя теперь теорему Ролля и учитывая \eqref{10} получаем, что
функция $\Delta_N^{(r-2)}(t)$ имеет не менее $2N+2$ нуля на
промежутке
$\left[\tau_1-N\cdot\frac{a+2}{\lambda},\tau_1+N\cdot\frac{a+2}{\lambda}\right]$.
Отсюда следует, что на некотором промежутке  монотонности
$[\alpha,\alpha+\frac{2+a}\lambda]\subset
\left[\tau_1-N\cdot\frac{a+2}{\lambda},
\tau_1+N\cdot\frac{a+2}{\lambda}\right]$ ($\alpha:= \frac
1\lambda+k\cdot\frac{a+2}\lambda,\,k\in \NN$) функции
$\Psi^{(r-2)}(t)=b\lambda^{-2}\psi (a, \lambda t)$ функция
$\Delta_N^{(r-2)}(t)$ меняет знак не менее трёх раз. Но тогда и
разность
$$
b\lambda^{-2}\psi (0, \lambda t)-x^{(r-2)}_N(t)
$$
на некотором промежутке монотонности функции $b\lambda^{-2}\psi (0,
\lambda t)$ меняет знак не менее трех раз. Однако это противоречит
теореме сравнения Колмогорова (см., например,~\cite[Предложение
5.5.3]{korn2}) так как эйлеров сплайн $b\lambda^{-2}\psi (0, \lambda
t)$ является функцией сравнения для функции $x^{(r-2)}_N(t)$.

Как следствие из теоремы \ref{th::0} получаем следующую теорему,
которую можно рассматривать как неравенство типа Колмогорова,
оценивающее норму $k$ - й производной функции, через норму самой
функции и нормы ее $(r-2)$ -й и $r$ -й производных.

\begin{theorem}\label{th::1}
 Пусть $r\in\NN$ и $0=k_1<k_2<k_3=r-2$, $k_4=r$ и $x\in L_{\infty,\infty}^r(\RR)$.
 Пусть числа $a,b,\lambda>0$ таковы, что для функции $\Psi_{a,b,\lambda}(t)$ выполняются соотношения \eqref{4'}.
 Тогда $$\left\|x^{(k_2)}\right\|\leq\left\|\Psi^{(k_2)}_{a,b,\lambda}\right\|.$$
 \end{theorem}

\medskip

{\bf 4. Решение задачи Колмогорова для случая $0=k_1<k_2<k_3=r-2$,
$k_4=r$}.
\begin{theorem}
Пусть заданы натуральные числа $r\geq 4$,  $0=k_1<k_2<k_3=r-2$, $k_4=r$ и действительные числа $M_{k_1},M_{k_2},M_{k_3},M_{k_4}>0$. Существует функция $x\in L_{\infty,\infty}^r(\RR)$, такая, что
\begin{equation}\label{11}
\left\|x^{(k_i)}\right\|=M_{k_i},\,\,i=1,2,3,4
\end{equation} тогда и только тогда, когда выполняются неравенства
$$a)\,\,M_{r-2}\leq \frac{\|\varphi_{2}\|}{\|\varphi_{r-k}\|^{\frac{2}{r-k}}}M_{k}^{\frac 2{r-k}}M_r^{\frac {r-k-2}{r-k}},$$
$$b)\,\, M_0\geq \|\Psi_r(M_{k_2},M_{r-2},M_r)\|,\phantom{aaaa}$$ где функция $\Psi_r$ определена в замечании~$\ref{z::1}$.
\end{theorem}

Необходимость условия $a)$ следует из неравенства Колмогорова,
необходи\-мость условия $b)$ следует из теоремы~\ref{th::1}.

Для доказательства достаточности заметим, что в случае выполнения
усло\-вий $a)$ и $b)$ для функции
$$x(t):=\Psi_r(M_{k_2},M_{r-2},M_r;t)+M_0-\left\|\Psi_r(M_{k_2},M_{r-2},M_r)\right\|$$
выполняются равенства~\eqref{11}.

Теорема доказана.